\documentstyle[11pt]{article}
\begin{document}
\newcommand{\p}{\parallel }
\makeatletter \makeatother
\newtheorem{th}{Theorem}[section]
\newtheorem{lem}{Lemma}[section]
\newtheorem{de}{Definition}[section]
\newtheorem{rem}{Remark}[section]
\newtheorem{cor}{Corollary}[section]
\renewcommand{\theequation}{\thesection.\arabic {equation}}

\title{{\bf The Equivariant Noncommutative Atiyah-Patodi-Singer Index Theorem}
\thanks{Partially supported by  MOEC and the 973 project.}}
\author{ Yong Wang \\
{\scriptsize \it Nankai Institute of Mathematics Tianjin 300071, P.R.China;}\\
{\scriptsize \it wangy581@nenu.edu.cn }}

\date{}
\maketitle

\begin{abstract} In [Wu], the noncommutative
Atiyah-Patodi-Singer index theorem was proved. In this paper, we
extend this theorem to the equivariant case.\\


\end{abstract}

\section{Introduction}

\quad In [APS], Atiyah-Patodi-Singer proved their famous
Atiyah-Patodi-Singer index theorem for manifolds with boundary. In
[D], Donnelly extended this theorem to the equivariant case by
modifying the Atiyah-Patodi-Singer original method. In [Z], Zhang
got this equivariant Atiyah-Patodi-Singer index theorem by using a
direct geometric method$^{[{\rm LYZ}]}$.
 In [Wu], Wu proved the
Atiyah-Patodi-Singer index theorem in the framework of
noncommutative geometry. To do so, he introduced the total eta
invariant (called the higher eta invariant in [Wu]) which is the
generalization of the classical Atiyah-Patodi-Singer eta
invariants$^{[{\rm APS}]},$ then proved its regularity by using the
Getzler symbol calculus$^{[{\rm G1}]}$ as adopted in [BF] and
computed
 its radius of convergence. Subsequently, he proved the variation
formula of eta cochains, using which he got the noncommutative
Atiyah-Patodi-Singer index theorem. In [G2], using superconnection,
Getzler gave another proof of the noncommutative
Atiyah-Patodi-Singer index theorem, which was more difficult, but
avoided mention of the operators $b$ and $B$ of
cyclic cohomology.\\
\indent The purpose of this paper is to extend the noncommutative
Atiyah-Patodi-Singer index theorem to the equivariant case.\\
\indent The paper is organized as follows: In Section 2.1, we define
the equivariant eta cochains and prove their regularity at infinity.
In Section 2.2
 we decompose the equivariant eta cochains
 into two parts and estimate the
second part. The first part will be estimated in Section 2.3. In
Section 2.4, we consider the convergence of the total eta
invariant.\\
 \indent Let $N$ be an odd dimensional spin manifold and $G$ be a compact Lie group acting
on $N$ by oritention-preserving isometries. let
$$C^{\infty}_G(N)=\{f\in C^{\infty}(N)|f(g.x)=f(x),~{\rm for~
any}~g\in G~ {\rm and}~ x\in N\}.$$
 Suppose that Dirac operator $D$ is invertible
with $\lambda$ the smallest eigenvalue of $|D|$, and $p=p^*=p^2\in
{\cal M }_r(C^{\infty}_G(N))$ is an idempotent which satisfying
$||dp||<\lambda$. Let $\eta^G(p(D\otimes
 I_r)p)$ is the equivariant Atiyah-Patodi-Singer eta
invariant associated to $p(D\otimes
 I_r)p$ which is the
Dirac operator with coefficients from $F=p({\bf C^r})$. $\eta^G(D)$
is the equivariant total eta invariants defined in Section 2.
Ch$(p)$ is the Chern character of $p$ defined in [GS].
 In Section 3, using the superconnction method in [G2], we
 prove the formula
 $$\frac{1}{2}\eta^G(p(D\otimes
 I_r)p)=\langle\eta^G(D),{\rm Ch}(p)\rangle,\eqno(1.1)$$
\indent In section 4, we define the equivariant Chern-Connes
character on manifolds with boundary and discuss its radius of
convergence.\\
\indent In section 5, we prove our main results. Using (1.1), we
express the equivariant index of the Dirac operator with the
coefficient from $G$-vector bundle $p({\bf C^r})$ over the cone as a
pair of the equivariant Chern-Connes character and ${\rm
Ch}(p)$.\\

\section{The Equivariant Total Eta Invariants }

{\bf 2.1 The Equivariant Eta Cochains}\\

 \quad Let $N$ be a compact
oriented odd dimensional Riemannian manifold without boundary with a
fixed spin structure and $S$ be the bundle of spinors on $N$. Denote
by $D$ the associated Dirac operator on $H=L^2(N;S)$, the Hilbert
space of $L^2$-sections of the bundle $S$. Let $c(df): S\rightarrow
S$ denote the Clifford action with $f\in C^{\infty}(N)$. Suppose
that $G$ is a connected Lie group acting on $N$ by
orientation-preserving isometries and $g\in G$ has a lift
$\widetilde{dg}:~\Gamma(S)\rightarrow\Gamma(S)$ (see [LYZ]), then we
have $\widetilde{dg}$ commutes with the Dirac operator and
$\widetilde{dg}$ is a bounded operator.

 Let
$A=C^{\infty}(N)$, then the data $(A,H,D,G)$ defines a finitely
(hence $\theta$-summable) equivariant unbounded Fredholm module in
the sense of [KL] (for details see [CH],[FGV] and [KL]). Similar to
[CH] or [W], for equivariant $\theta$-summable Fredholm module
$(A,H,D,G)$, we can define the equivariant cochain $\widetilde{{\bf
{\rm ch}}_
k^G}(tD,D)$ ($k$ is even) by the formula:\\

$~~~~~~~~~~\widetilde{{\bf {\rm ch}}_k^G}(tD,D)(f^0,\cdots ,
f^k)(g)$
$$:=
t^k\sum^k_{i=0}(-1)^i\langle f^0,c(df^1),\cdots ,
c(df^i),D,c(df^{i+1}),\cdots , c(df^k)\rangle _t(g),\eqno(2.1.1)$$
where $f^0,\cdots,f^k \in C^{\infty}(N),~ g\in G$. If $A_i~ (0\leq i\leq n)$ are operators on $H$, we define:\\
\bigskip
$$\langle A_0,\cdots , A_n\rangle_t(g)=\int_{\triangle_n}{\rm Tr}(A_0e^{-t^2s_1D^2}A_1e^{-t^2(s_2-s_1)D^2}\cdots A_n
e^{-t^2(1-s_n)D^2}\widetilde{dg})ds,\eqno(2.1.2)$$ where
$\triangle_n=\{(s_1,\cdots,s_n)|~0\leq
s_1\leq\cdots\leq s_n\leq 1\}$ is the simplex in ${\bf R^n}$.\\
\indent Formally, the  equivariant total $\eta$-invariant of the
Dirac operator $D$ is defined to be a sequence of even equivariant
cochains on $C^{\infty}(N)$,
by the formula:\\
$$\eta_k^G(D)=\frac{1}{\Gamma(\frac{1}{2})}\int^{\infty}_{0}\widetilde{{\bf {\rm ch}}_k^G}(tD,D)dt,
\eqno(2.1.3)$$ \noindent where $\Gamma(\frac{1}{2})=\sqrt {\pi}.$
Then $\eta_0^G(D)(1)(g)$ is the half of the equivariant eta
invariants defined in [APS], [D] and [Z]. In order to prove that the
above definition is well defined, it is necessary to check the
integrality near the two ends of the integration.
 Firstly, the regularity at infinity comes from the
following lemma.\\

 \noindent {\bf Lemma 2.1}~~{\it For
$f^0,\cdots,f^k \in C^{\infty}(N)$ and $g\in G$, we have}
$$\widetilde{{\bf {\rm ch}}_k^G}(tD,D)(f^0,\cdots , f^k)(g)=O(t^{-2}),~~{\rm as}~t\rightarrow\infty.\eqno(2.1.4)$$
\noindent{\bf {\it Proof.}}~~Since $\widetilde{dg}$ is a bounded
operator, our proof is similar to the proof of Lemma 2 in [CH]. ~~~$\Box$\\

\noindent{\bf 2.2 Expansion of The Equivariant  Eta Cochains }\\

\quad In [W], Wu proved the regularity at zero of (2.1.3) in the
$g=id$ case by using the Getzler symbol calculus. In what follows,
we will give a proof of the regularity of (2.1.3) at zero in the
general case by using the method in [CH] and [F].

Firstly, recall some Lemmas in [CH] and [F].

Let $H$ be a Hilbert space. For $q\geq0$, denote by $||.||_q$
Schatten $p$-norm on Schatten ideal $L^p$ (for details, see [S]).
$L(H)$ denotes the Banach algebra of bounded operators on $H$.\\

\noindent {\bf Lemma 2.2}~([CH],[F]){\it~~(i)~~${\rm Tr}(AB)={\rm
Tr}(BA)$, for $A,~B\in L(H)$ and $AB, ~BA\in
L^1$.\\
~~(ii)~~For $A\in L^1,$ we have $|{\rm Tr}(A)|\leq ||A||_1$,
$||A||\leq ||A||_1$.\\
~~(iii)~~For $A\in L^q$ and $B\in L(H)$, we have: $||AB||_q\leq
||B||||A||_q$, $||BA||_q\leq ||B||||A||_q$.\\
 ~~(iv)~(H\"{o}lder
Inequality)~~If $\frac{1}{r}=\frac{1}{p}+\frac{1}{q},~p,q,r>0,~A\in
L^p,~B\in
L^q,$ then $AB\in L^r$ and $||AB||_r\leq ||A||_p||B||_q$.}\\

\noindent{\bf Lemma 2.3}~([CH],[F])~{\it  For any $u>0,~t>0$ and any
order $l$
differential operator $B$, we have:}\\
$$||e^{-utD^2}B||_{u^{-1}}\leq C_lu^{-\frac{l}{2}}t^{-\frac{l}{2}}({\rm
tr}[e^{-\frac{tD^2}{2}}])^u.\eqno(2.2.1)$$\\

 \noindent{\bf Lemma 2.4}~([CH],[F]) {\it ~Let $B_1,~B_2$ be positive order
 $p,~q$ pseudodifferential operators respectively, then for any
 $s,~t>0,~0\leq u\leq1$, we have the following estimate:}\\
$$||B_1e^{-ustD^2}B_2e^{-(1-u)stD^2}||_{s^{-1}}\leq C_{p,q}s^{-\frac{p+q}{2}}t^{-\frac{p+q}{2}}({\rm
tr}[e^{-\frac{tD^2}{4}}])^s.\eqno(2.2.2)$$
 \indent Let $B$ be an
operator and $l$ be a positive interger.
 Write
$$B^{[l]}=[D^2,B^{[l-1]}],~B^{[0]}=B.$$\\
\noindent{\bf Lemma 2.5}~([CH],[F])~{\it Let $B$ a finite order
differential
operator, then for any $s>0$, we have:}\\
$$e^{-sD^2}B=\sum^{N-1}_{l=0}\frac{(-1)^l}{l!}s^lB^{[l]}e^{-sD^2}+(-1)^Ns^NB^{[N]}(s),\eqno(2.2.3)$$
\noindent {\it where $B^{[N]}(s)$ is given by}\\
$$B^{[N]}(s)=\int_{\triangle_N}e^{-u_1sD^2}B^{[N]}e^{-(1-u_1)sD^2}du_1du_2\cdots du_N.\eqno(2.2.4)$$
\indent Similar to Lemma 5, we have:\\

\noindent{\bf Lemma 2.6}~ {\it ~~Let $B$ a finite order differential
operator, then for any $s>0$, we have:}\\
$$Be^{-sD^2}=\sum^{N-1}_{l=0}\frac{1}{l!}s^le^{-sD^2}B^{[l]}+s^NB_1^{[N]}(s),\eqno(2.2.5)$$
\noindent {\it where $B_1^{[N]}(s)$ is given by}\\
$$B_1^{[N]}(s)=\int_{\triangle_N}e^{-(1-u_1)sD^2}B^{[N]}e^{-u_1sD^2}du_1du_2\cdots du_N.\eqno(2.2.6)$$
\indent In order to prove that (2.1.3) is well defined, it is enough
to prove that when $t\rightarrow 0$, we have the estimate
$$\widetilde{{\bf {\rm ch}}_k^G}(tD,D)(f^0,\cdots ,f^k)(g)\sim O(t);~~ i.e.~~
\widetilde{{\bf {\rm ch}}_k^G}(\sqrt {t}D,D)(f^0,\cdots ,f^k)(g)\sim
O(t^{\frac{1}{2}}).\eqno(2.2.7)$$
 By (2.1.5), we have for
$i=0,\cdots ,k$, the $i-$th term of
$\widetilde{{\bf {\rm ch}}_k^G}(\sqrt {t}D,D)(f^0,\cdots ,f^k)(g)$ up to sign is:\\
$$t^{\frac{k}{2}}\langle f^0,c(df^1),\cdots ,c(df^i),D,c(df^{i+1}),\cdots ,c(df^k)\rangle _{\sqrt{t}}(g)$$
$$=
t^{\frac{k}{2}}\int_{\triangle_{k}}{\rm
Tr}[f^0e^{-s_1tD^2}c(df^1)e^{-(s_2-s_1)tD^2}c(df^2)\cdots
c(df^i)(s_{i+1}-s_i) De^{-(s_{i+1}-s_i)tD^2}$$ $$\cdot
c(df^{i+1})\cdots c(df^k) e^{-(1-s_k)tD^2}\widetilde{dg}]ds_1\cdots
ds_k.\eqno(2.2.8)$$ \quad In what follows, we will compute the
expression of (2.2.8) by the above lemmas. By
Lemma 2.6, we have:\\
$$c(df^{i+1})e^{-(s_{i+2}-s_{i+1})tD^2}\cdots c(df^k)
e^{-(1-s_k)tD^2}$$
$$=\sum^{N-1}_{\lambda _{i+1},\cdots ,\lambda _k=0}\frac{(1-s_{i+1})^{\lambda _{i+1}}\cdots (1-s_k)^{\lambda
_k}t^{\lambda _{i+1}+\cdots +\lambda _k}}{\lambda
_{i+1}!\cdots\lambda _k!}$$  $$\cdot
e^{-(1-s_{i+1})tD^2}[c(df^{i+1})]^{[\lambda _{i+1}]}\cdots
[c(df^k)]^{[\lambda _k]}$$
$$+\sum_{i+1\leq q\leq k}\sum^{N-1}_{\lambda _{q+1},\cdots ,\lambda _k=0}
\frac{(1-s_q)^N(1-s_{q+1})^{\lambda _{q+1}}\cdots (1-s_k)^{\lambda
_k} t^{N+\lambda _{q+1}+\cdots +\lambda _k}}{\lambda
_{q+1}!\cdots\lambda _k!}$$
$$\cdot c(df^{i+1})e^{-(s_{i+2}-s_{i+1})tD^2}\cdots c(df^{q-1})
e^{-(s_q-s_{q-1})tD^2}$$ $$\cdot
\{[c(df^q)]_1^{N}((1-s_q)t)\}[c(df^{q+1})]^{[\lambda_{q+1}]}\cdots
[c(df^k)]^{[\lambda_k]}.\eqno(2.2.9)$$
\noindent By Lemma 2.5, we have:\\
$$f^0e^{-s_1tD^2}c(df^1)e^{-(s_2-s_1)tD^2}c(df^2)\cdots c(df^i)
e^{-(s_{i+1}-s_i)tD^2}$$
$$=\sum^{N-1}_{\lambda _{1},\cdots ,\lambda _i=0}\frac{(-1)^{\lambda _1+\cdots +\lambda _i}{s_1}^{\lambda _{1}}\cdots
 s_i^{\lambda
_i}t^{\lambda _{1}+\cdots +\lambda _i}}{\lambda _{1}!\cdots\lambda
_i!}f^0[c(df^{1})]^{[\lambda _{1}]}\cdots  [c(df^i)]^{[\lambda
_i]}e^{-s_{i+1}tD^2}$$
$$+\sum_{1\leq q\leq i}\sum^{N-1}_{\lambda _{1},\cdots ,\lambda _{q-1}=0}
\frac{(-1)^{\lambda _1+\cdots +\lambda _{q-1}+N}s_1^{\lambda
_1}\cdots s_{q-1}^{\lambda _{q-1}}s_q^Nt^{\lambda _1+\cdots +\lambda
_{q-1}+N}}{\lambda _{1}!\cdots\lambda
_{q-1}!}f^0[c(df^{1})]^{[\lambda _{1}]}$$
$$\cdots[c(df^{q-1})]^{[\lambda _{q-1}]}\{ [c(df^q)]^{[N]}(s_qt)\}
e^{-(s_{q+1}-s_q)tD^2}\cdots
c(df^{i})e^{-(s_{i+1}-s_i)tD^2}.\eqno(2.2.10)$$
\noindent By (2.2.8) and (2.2.10), we have:\\
$$t^{\frac{k}{2}}(s_{i+1}-s_i){\rm Tr}[f^0e^{-s_1tD^2}c(df^1)e^{-(s_2-s_1)tD^2}c(df^2)\cdots c(df^i)
De^{-(s_{i+1}-s_i)tD^2}$$ $$c(df^{i+1})\cdots c(df^k)
e^{-(1-s_k)tD^2}\widetilde{dg}]$$
$$=\sum^{N-1}_{\lambda _{1},\cdots ,\lambda _i=0}\frac{(-1)^{\lambda _1+\cdots +\lambda _i}{s_1}^{\lambda _{1}}\cdots
 s_i^{\lambda
_i}(s_{i+1}-s_i)t^{\lambda _{1}+\cdots +\lambda
_i+\frac{k}{2}}}{\lambda _{1}!\cdots\lambda _i!}{\rm
Tr}\{f^0[c(df^{1})]^{[\lambda _{1}]}$$ $$\cdots  [c(df^i)]^{[\lambda
_i]}e^{-s_{i+1}tD^2}Dc(df^{i+1})\cdots c(df^k)
e^{-(1-s_k)tD^2}\widetilde{dg}\}$$
$$+\sum_{1\leq q\leq i}\sum^{N-1}_{\lambda _{1},\cdots ,\lambda _{q-1}=0}
\frac{(-1)^{\lambda _1+\cdots +\lambda _{q-1}+N}s_1^{\lambda
_1}\cdots s_{q-1}^{\lambda _{q-1}}s_q^{N}(s_{i+1}-s_i) t^{\lambda
_1+\cdots +\lambda _{q-1}+N+\frac{k}{2}}}{\lambda _{1}!\cdots\lambda
_{q-1}!}$$
$$\cdot {\rm Tr}\{f^0[c(df^{1})]^{[\lambda _{1}]}\cdots[c(df^{q-1})]^{[\lambda _{q-1}]}\{
[c(df^q)]^{[N]}(s_qt)\} e^{-(s_{q+1}-s_q)tD^2}\cdots$$
$$\cdot c(df^{i})e^{-(s_{i+1}-s_i)tD^2}Dc(df^{i+1})\cdots c(df^k)
e^{-(1-s_k)tD^2}\widetilde{dg}\}$$
$$=\sum_{0\leq \lambda _{1},\cdots ,\lambda _k\leq {N-1}}\frac{(-1)^{|\lambda|}{s_1}^{\lambda _{1}}\cdots s_i^{\lambda
_i}(s_{i+1}-1)^{\lambda _{i+1}}\cdots (s_k-1)^{\lambda
_k}t^{|\lambda|+\frac{k}{2}}(s_{i+1}-s_i)}{\lambda!}$$
$$\cdot {\rm Tr}\{f^0[c(df^{1})]^{[\lambda _{1}]}\cdots [c(df^i)]^{[\lambda
_i]}De^{-tD^2}[c(df^{i+1})]^{[\lambda _{i+1}]}\cdots
[c(df^k)]^{[\lambda _k]}\widetilde{dg}\}+A^i_1+A^i_2,\eqno(2.2.11)$$
\noindent where the last equality comes from (2.2.9) and\\
$$A^i_1=\sum_{1\leq q\leq i}\sum^{N-1}_{\lambda _{1},\cdots ,\lambda _{q-1}=0}
\frac{(-1)^{\lambda _1+\cdots +\lambda
_{q-1}+N}(s_{i+1}-s_i)s_1^{\lambda _1}\cdots s_{q-1}^{\lambda
_{q-1}}s_q^{N}t^{\lambda _1+\cdots +\lambda _{q-1}+N+\frac{k}{2}}}
{\lambda _{1}!\cdots\lambda _{q-1}!}$$
$$\cdot {\rm Tr} \{f^0 [c(df^{1})]^{[\lambda
_{1}]}\cdots [c(df^{q-1})]^{[\lambda _{q-1}]}\{
[c(df^q)]^{[N]}(s_qt)\} e^{-(s_{q+1}-s_q)tD^2}\cdots $$ $$\cdot
c(df^{i})e^{-(s_{i+1}-s_i)tD^2}D c(df^{i+1})\cdots c(df^k)
e^{-(1-s_k)tD^2}\widetilde{dg}\};\eqno(2.2.12)$$
\\
$~~A^i_2=\sum\limits^{N-1}_{\lambda _{1},\cdots ,\lambda _i=0}
\sum\limits_{i+1\leq q\leq k}\sum\limits^{N-1}_{\lambda
_{q+1},\cdots ,\lambda _k=0}\frac{(-1)^{\lambda _1+\cdots +\lambda
_i}}{\lambda _{1}!\cdots\lambda _i!\lambda _{q+1}!\cdots\lambda
_k!} (s_{i+1}-s_i){s_1}^{\lambda _{1}}\cdots s_i^{\lambda _i}$\\
$$ \cdot(1-s_q)^{N}(1-s_{q+1})^{\lambda _{q+1}}\cdots
(1-s_k)^{\lambda _k}  t^{\lambda _{1}+\cdots +\lambda _i+N+\lambda
_{q+1}+\cdots+\lambda _k+\frac{k}{2}}$$ $$\cdot{\rm Tr}
\{f^0[c(df^{1})]^{[\lambda _{1}]}\cdots [c(df^i)]^{[\lambda
_i]}e^{-s_{i+1}tD^2}D c(df^{i+1}) e^{-(s_{i+2}-s_{i+1})tD^2}\cdots
c(df^{q-1})$$ $$\cdot e^{-(s_q-s_{q-1})tD^2}\{
[c(df^q)]_1^{[N]}((1-s_q)t)\} [c(df^{q+1})]^{[\lambda_{q+1}]}\cdots
[c(df^k)]^{[\lambda_k]}\widetilde{dg}\}.\eqno(2.2.13)$$ \noindent
Using Lemma 2.2 (i) by taking appropriate bounded operators $A,B$,
we have: $${\rm Tr}\{f^0[c(df^{1})]^{[\lambda _{1}]}\cdots
[c(df^i)]^{[\lambda _i]}De^{-tD^2}[c(df^{i+1})]^{[\lambda
_{i+1}]}\cdots [c(df^k)]^{[\lambda _k]}\widetilde{dg}\}$$
\begin{eqnarray*}
&=&{\rm Tr}\{[c(df^{i+1})]^{[\lambda _{i+1}]}\cdots
[c(df^k)]^{[\lambda _k]}\widetilde{dg}f^0[c(df^{1})]^{[\lambda
_{1}]}\cdots [c(df^i)]^{[\lambda _i]}De^{-tD^2}\}\\
&=&{\rm Tr}\{D^{\lambda}_iDe^{-tD^2}\},
~~~~~~~~~~~~~~~~~~~~~~~~~~~~~~~~~~~~~~~~~~~~~~~~~~~~~~~~~~~~~~~~~~~~~~(2.2.14)
\end{eqnarray*}
where we use the notation\\
$$D^{\lambda}_i=[c(df^{i+1})]^{[\lambda
_{i+1}]}\cdots [c(df^k)]^{[\lambda
_k]}\widetilde{dg}f^0[c(df^{1})]^{[\lambda _{1}]}\cdots
[c(df^i)]^{[\lambda _i]}.\eqno(2.2.15)$$ Using
(2.1.1),(2.2.11),(2.2.12),(3.13) and (2.2.14), then we have:\\

\noindent{\bf Corollary 2.7}~~{\it Let {\rm dim}$M=n=2m+1$, for
$f_j\in
C^{\infty}(M), ~0\leq j\leq k,~k$ is even and $g\in G$, then}\\
$\widetilde{{\bf {\rm ch}}_k^G}(\sqrt {t}D,D)(f^0,\cdots ,f^k)(g)$\\
$$=\sum^k_{i=0}(-1)^i \sum_{0\leq \lambda _{1},\cdots ,\lambda _k\leq
{N-1}}\frac{(-1)^{|\lambda|}Ct^{|\lambda|+\frac{k}{2}}}{\lambda!}
{\rm
Tr}\{D^{\lambda}_iDe^{-tD^2}\}+\sum^k_{i=0}(-1)^i\int_{\triangle_k}(A_1^i+A_2^i)ds,\eqno(2.2.16)$$
{\it with the constant}\\
$$C={\sum\limits^{\lambda_{i+1}}_{j_{i+1}=0}}\cdots\sum\limits^{\lambda_k}_{j_k=0}(\prod\limits^k_{s=i+1}
\left(\begin{array}{lcr}
  \ \lambda_s  \\
    \  j_s
\end{array}\right)
(-1)^{\lambda_s-j_s}) \frac{1}{\lambda_1+1}\cdots
\frac{1}{\sum\limits^{i-1}_{j=1}{\lambda_j}+i-1}
\frac{1}{\sum\limits^{i}_{j=1}{\lambda_j}+i}$$
$$\cdot\frac{1}{\sum\limits^{i}_{j=1}{\lambda_j}+i+1}
\frac{1}{\sum\limits^{i}_{j=1}{\lambda_j}+j_{i+1}+i+2}\cdots
\frac{1}{\sum\limits^{i}_{j=1}{\lambda_j}+j_{i+1} +\cdots
+j_k+k+1},$$ \indent Now we give the estimate of $\widetilde{{\bf
{\rm ch}}_k^G}(\sqrt {t}D,D)$ and
let $N=n+2-k=2m+3-k$ in (2.2.16). We have:\\

\noindent{\bf Theorem 2.8}~~{\it 1)~~If $k\leq$ {\rm dim}$N+1=2m+2$,
then
when $t\rightarrow 0^+$, we have:}\\
$~~~~~~~~~~~~~~~\widetilde{{\bf {\rm ch}}_k^G}(\sqrt
{t}D,D)(f^0,\cdots ,f^k)(g)$
$$=\sum^k_{i=0}(-1)^i
\sum_{0\leq \lambda _{1},\cdots ,\lambda _k\leq
{N-1}}\frac{(-1)^{|\lambda|}Ct^{|\lambda|+\frac{k}{2}}}{\lambda!}
{\rm
Tr}\{D^{\lambda}_iDe^{-tD^2}\}+O(t^{\frac{1}{2}}).\eqno(2.2.17)$$
{\it 2)~~If $k>2m+2$, then
when $t\rightarrow 0^+$, we have:}\\
$$\widetilde{{\bf {\rm ch}}_k^G}(\sqrt {t}D,D)(f^0,\cdots ,f^k)(g)\sim O(t^{\frac{1}{2}}).\eqno(2.2.18)$$
\noindent {\bf {\it Proof.}}~~1)~In order to prove (2.2.17), we
 only prove that when $t\rightarrow 0^+$,~$\int_{\triangle_k}A_2^ids\sim
O(t^{\frac{1}{2}})$ (similar $\int_{\triangle_k}A_1^ids\sim
O(t^{\frac{1}{2}})$). By (2.2.13), then
\begin{eqnarray*}
|\int_{\triangle_k}A_2^ids|
&\leq&\int_{\triangle_k}|A_2^i|ds\\
&\leq&\sum^{N-1}_{\lambda _{1},\cdots ,\lambda _i=0} \sum_{i+1\leq
q\leq k}\sum^{N-1}_{\lambda _{q+1},\cdots ,\lambda
_k=0}\frac{A_{\lambda _1,\cdots ,\lambda _i,\lambda
_{q+1},\cdots,\lambda _k}}{\lambda _{1}!\cdots\lambda _i!\lambda
_{q+1}!\cdots\lambda _k!},~~~~~~~~(2.2.19)
\end{eqnarray*}
where\\
\begin{eqnarray*}
A_{\lambda _1,\ldots,\lambda _i,\lambda _{q+1},\ldots,\lambda_k}
&=&\int_{\triangle_k}{(s_{i+1}-s_i){s_1}^{\lambda_{1}}\cdots
s_i^{\lambda _i}(1-s_q)^{N}(1-s_{q+1})^{\lambda _{q+1}}
\cdots (1-s_k)^{\lambda_k}}\\
&&\cdot t^{\lambda _{1}+\cdots +\lambda _i+N+\lambda _{q+1}+\cdots
+\lambda_k+\frac{k}{2}} \left|{\rm Tr} \{f^0[c(df^{1})]^{[\lambda
_{1}]}\cdots [c(df^i)]^{[\lambda _i]}\right.\\
&&\cdot e^{-s_{i+1}tD^2}D c(df^{i+1}) e^{-(s_{i+2}-s_{i+1})tD^2}
\cdots c(df^{q-1})e^{-(s_q-s_{q-1})tD^2}\\
&& \{ [c(df^q)]_1^{[N]}((1-s_q)t)\} [c(df^{q+1})]^{[\lambda_{q+1}]}
\end{eqnarray*}
$$\left.\cdots [c(df^k)]^{[\lambda_k]}\widetilde{dg}\}\right|ds.
\eqno(2.2.20)$$
\noindent By Lemma 2.2, Lemma 2.4 and Lemma 2.6, then:\\
$$A_{\lambda _1,\ldots,\lambda _i,\lambda _{q+1},\ldots,\lambda_k}\leq
\int_{\triangle_k}{(s_{i+1}-s_i){s_i}^{\lambda _{1}+\cdots +\lambda
_i}(1-s_q)^{N}(1-s_{q+1})^{\lambda _{q+1}} \cdots
(1-s_k)^{\lambda_k}}$$ $$\cdot t^{\lambda _{1}+\cdots +\lambda
_i+N+\lambda _{q+1}+\cdots +\lambda_k+\frac{k}{2}} \left|{\rm
Tr}\left\{
 \left(f^0[c(df^{1})]^{[\lambda _{1}]}\cdots [c(df^i)]^{[\lambda
_i]}e^{-s_{i}tD^2}\right)
\left(De^{-(s_{i+1}-s_{i})tD^2}\right)\right.\right.$$
$$.\left(
c(df^{i+1}) e^{-(s_{i+2}-s_{i+1})tD^2}\right)\cdots\left(c(df^{q-1})
e^{-(s_q-s_{q-1})tD^2}\right)$$
$$\cdot\left.\left.\left(\{
[c(df^q)]_1^{[N]}((1-s_q)t)\} [c(df^{q+1})]^{[\lambda _{q+1}]}\cdots
[c(df^k)]^{[\lambda_k]}\widetilde{dg}\right) \right\}\right|ds$$
$$\leq C_0\int_{\triangle_q}\int_{\triangle_N}(s_{i+1}-s_{i}){s_i}^{\lambda _{1}+\cdots +\lambda
_i}(1-s_q)^{N+\lambda _{q+1}+\cdots\lambda _k+k-q}t^{N+\lambda
_{1}+\cdots +\lambda _i+\lambda _{q+1}+\cdots +\lambda
_k+\frac{k}{2}}$$
$$\cdot
\left|\left| \left(f^0[c(df^{1})]^{[\lambda _{1}]}\cdots
[c(df^i)]^{[\lambda _i]}e^{-s_{i}tD^2}\right)
\left(De^{-(s_{i+1}-s_{i})tD^2}\right)\left( c(df^{i+1})
e^{-(s_{i+2}-s_{i+1})tD^2}\right)\right.\right.$$
$$\cdots\left(c(df^{q-1})
e^{-(s_q-s_{q-1})tD^2}\right) \cdot
\left(e^{-(1-u_1)(1-s_q)tD^2}[c(df^q)]^{[N]}e^{-u_1(1-s_q)tD^2}\right.$$
$$\cdot\left.\left.\left.[c(df^{q+1})]^{[\lambda_{q+1}]}\cdots
[c(df^k)]^{[\lambda_k]}\widetilde{dg}\right)
\right|\right|_1du_1\cdots du_Nds_1\cdots ds_q$$
$$\leq C_0\int_{\triangle_q}\int_{\triangle_N}
t^{N+\lambda _{1}+\cdots +\lambda _i+\lambda _{q+1}+\cdots +\lambda
_k+\frac{k}{2}} \left|\left|{s_i}^{\lambda _{1}+\cdots +\lambda
_i}f^0[c(df^{1})]^{[\lambda _{1}]}\cdots [c(df^i)]^{[\lambda
_i]}e^{-s_{i}tD^2}\right|\right|_{\frac{1}{s_i}}$$
$$\cdot\left|\left|(s_{i+1}-s_{i})De^{-\frac{(s_{i+1}-s_{i})}{2}tD^2}\right|\right|
\left|\left|e^{-\frac{(s_{i+1}-s_{i})}{2}tD^2}\right|\right|_{\frac{1}{s_{i+1}-s_i}}
\left|\left|c(df^{i+1}) e^{-(s_{i+2}-s_{i+1})tD^2}
\right|\right|_{\frac{1}{s_{i+2}-s_{i+1}}}$$
$$\cdots
 \left|\left|c(df^{q-1})
e^{-(s_q-s_{q-1})tD^2}\right|\right|_{\frac{1}{s_q-s_{q-1}}}
\left|\left|(1-s_q)^{N+\lambda _{q+1}+\cdots\lambda
_k+k-q}e^{-(1-u_1)(1-s_q)tD^2}\right.\right.$$
$$\cdot\left.\left.[c(df^q)]^{[N]}e^{-u_1(1-s_q)tD^2}
[c(df^{q+1})]^{[\lambda_{q+1}]}\cdots
[c(df^k)]^{[\lambda_k]}\widetilde{dg}\right|\right|_{\frac{1}{1-s_q}}$$
$$\leq \overline{C}\int_{\triangle_q}\int_{\triangle_N}
t^{\frac{N+\lambda _{1}+\cdots +\lambda _i+\lambda _{q+1}+\cdots
+\lambda _k+k-1}{2}}{s_i}^{\frac{\lambda _{1}+\cdots +\lambda
_i}{2}}\sqrt {s_{i+1}-s_{i}}(1-s_q)^{\frac{N+\lambda
_{q+1}+\cdots\lambda _k}{2}+k-q}$$ \noindent$$.{\rm
Tr}\{e^{-\frac{tD^2}{4}}\}du_1\cdots du_Nds_1\cdots ds_q,\eqno(2.2.21)$$\\
\noindent where $C_0,~\overline{C}$ are constants and the second
inequality comes from integrating respect to $s_{q+1},\cdots , s_k$
and (2.2.6). In the third inequality we use Lemma 2.2 (iv) and the
last inequality comes from Lemma 2.2 and Lemma 2.4. By the Weyl
asymptotics on the heat kernel we have when $t\rightarrow 0$,
$${\rm
Tr}\{e^{-\frac{tD^2}{4}}\}\sim O(t^{-\frac{n}{2}}).\eqno(2.2.22)$$
By (2.2.22) and $N=n+2-k$, So (2.2.21) $\sim O(t^{\frac{1}{2}}).$

 2)~~If $k>{\rm dim}N+1$, then
$$\left|t^{\frac{k}{2}}\int_{\triangle_k}(s_{i+1}-s_i){\rm Tr}[f^0e^{-s_1tD^2}c(df^1)e^{-(s_2-s_1)tD^2}c(df^2)\cdots
c(df^i) \right.$$
$$\left. \cdot De^{-(s_{i+1}-s_i)tD^2} c(df^{i+1})\cdots c(df^k)
e^{-(1-s_k)tD^2}\widetilde{dg}]ds\right|$$
$$\leq
t^{\frac{k}{2}}\int_{\triangle_k}\left|\left|\left(f^0e^{-s_1tD^2}c(df^1)\right)
\left(e^{-(s_2-s_1)tD^2}c(df^2)\right)\cdots\left(e^{-(s_{i}-s_{i-1})tD^2}c(df^i)\right)\right.\right.$$
$$.\left((s_{i+1}-s_i)De^{-\frac{s_{i+1}-s_i}{2}tD^2}\right)
\left(e^{-\frac{s_{i+1}-s_i}{2}tD^2}c(df^{i+1})\right)\cdots$$
$$.\left.\left.\left((e^{-(s_{k}-s_{k-1})tD^2}c(df^k)\right)
\left(e^{-(1-s_k)tD^2}\right)\widetilde{dg}\right|\right|_1ds$$
$$\leq
t^{\frac{k}{2}}\int_{\triangle_k}\sqrt
{s_{i+1}-s_{i}}~t^{-\frac{1}{2}}{\rm
Tr}\{e^{-\frac{tD^2}{4}}\}ds~\sim O(t^{\frac{-n+k-1}{2}})\sim
O(t^{\frac{1}{2}}),$$ \noindent where we use Lemma 2.2, Weyl
estimate and condition $k>{\rm
dim}N+1$.~~~$\Box$\\

\noindent {\bf 2.3 Clifford Asymptotics for Heat Kernels}\\

\indent By Theorem 2.8, in order to prove the regularity of the
equivariant eta cochains, it is enough to prove that when
$t\rightarrow 0$,
$$t^{|\lambda|+\frac{k}{2}}{\rm
Tr}\{D^{\lambda}_iDe^{-tD^2}\}\sim O(t^{\frac{1}{2}}).~~~~~~$$
\noindent Similar to Theorem 1.1 in [Z], we have the following lemma.\\

\noindent{\bf Lemma 2.9}~~$${\rm
Tr}\{D^{\lambda}_iDe^{-tD^2}\}=\int_M{\rm
Tr}\{(D^{\lambda}_i)_x[D{\rm{exp}}(-tD^2)(x,y)]|_{y=g\cdot
x}\}dx.\eqno(2.3.1)$$ \\

\noindent{\bf Proposition 2.10}~~{\it If $g$ has no fixed points on
$N$, then
$${\rm lim}_{t\rightarrow 0}t^{-\frac{1}{2}}\int_M{\rm
Tr}\{(D^{\lambda}_i)_x[D{\rm{exp}}(-tD^2)(x,g\cdot x)]\}dx=0.
\eqno(2.3.2)$$} \noindent{\bf {\it Proof.}}~~We introduce an
auxiliary Grassmann variable $z$ as in [BF]. By Duhamel principle,
we have
$${\rm exp}(-t(D^2-zD))={\rm exp}(-tD^2)+ztD{\rm exp}(-tD^2).\eqno(2.3.3)$$
Since $g$ has no fixed points, $d(x,g\cdot x)>\delta$ for some
constant $\delta>0$. So there exist positive constants $C_i
~(i=1,2,3,4)$ and positive integers $m_1,m_2$ such that
$t\rightarrow 0$,
$$||(D^{\lambda}_i)_x{\rm{exp}}(-tD^2)(x,g\cdot x)||\leq \frac{C_1}
{t^{\frac{n}{2}+m_1}}\rm{exp}(-\frac{C_2}{t});\eqno(2.3.4)$$
$$||(D^{\lambda}_i)_x{\rm{exp}}(-t(D^2-zD))(x,g\cdot x)||\leq
\frac{C_3}{t^{\frac{n}{2}+m_2}}\rm{exp}(-\frac{C_4}{t}),\eqno(2.3.5)$$
then similar to Corollary 1.4 in [Z], we prove this Proposition.~~$\Box $\\
\indent Since $g$ is an isometry, the fixed point set $F$ of $g$
consists of components $F_1,\cdots,F_k$, each of even codimension.
If $U$ is an open neighborhood of $F$, then by Proposition 2.10, we
have
$${\rm lim}_{t\rightarrow 0}t^{-\frac{1}{2}}\int_M{\rm
Tr}\{(D^{\lambda}_i)_x[D{\rm{exp}}(-tD^2)(x,g\cdot x)]\}dx$$
$$={\rm lim}_{t\rightarrow 0}t^{-\frac{1}{2}}\int_U{\rm
Tr}\{(D^{\lambda}_i)_x[D{\rm{exp}}(-tD^2)(x,g\cdot x)]\}dx.
\eqno(2.3.6)$$ We may assume $k=1$ and codim$F=2n'$. Denote by
$N(F)$ the normal bundle to $F$, similar to Theorem 2.2 in [LYZ] we
need only to prove that
$${\rm lim}_{t\rightarrow 0}t^{-\frac{1}{2}}|\int_F\int_{N_\xi(\varepsilon)}t^{|\lambda|+\frac{k}{2}}{\rm
Tr}\{(D^{\lambda}_i)_x[D{\rm{exp}}(-tD^2)(x,g\cdot x)]\}dN_\xi
d\xi|\leq C,\eqno(2.3.7)$$
\noindent for some constant $C>0$. Here $N_\xi(\varepsilon)=\{v\in N_\xi(F)|~||v||<\varepsilon~\}$.\\
\indent Similar to [LYZ],[Y], for $\xi\in F$, we choose an open
neighborhood $U$ of $\xi$ and the orthogonal frame over $U$ (see
[LYZ], pp.574). Consider the oriented
 orthonormal frame field $E^{g\cdot x}$ defined over the patch $U$ by
requiring that $E^{g\cdot x}(g\cdot x)=E(g\cdot x)$ and that
$E^{g\cdot x}$ is parallel along geodesics through
 $g\cdot x$. Choose a spin frame field $\sigma:U\rightarrow {\rm Spin}(M)$ such that $\pi\sigma=(E_1^{g\cdot x},\cdots,
E_n^{g\cdot x})$ where $\pi:{\rm Spin}(M)\rightarrow {\rm SO}(M)$ is
the two-fold covering over ${\rm SO}(M)$. For $x\in U$, let
$K(x),~g^*(x),~T^{\lambda}_i(x)\in {\rm Hom}(I,I)$ ($I$ is the
canonical spinors space as in [LYZ].) be defined through the
equivalence relations:
$$D{\rm{exp}}(-tD^2)(x,g\cdot x)[(\sigma(g\cdot x),u)]=[(\sigma(x),K(x)u)];\eqno(2.3.8)$$
$$\widetilde{dg}[(\sigma(x),v)]=[(\sigma(g\cdot x),g^*(x)v)];\eqno(2.3.9)$$
$$(D^{\lambda}_i)_x[(\sigma(x),w)]=[(\sigma(g\cdot x),T^{\lambda}_i(x)w)].\eqno(2.3.10)$$
Then similar to Lemma 4.1 in [LYZ], we have\\

 \noindent{\bf Lemma
2.11}~~$${\rm Tr}\{(D^{\lambda}_i)_x[D{\rm{exp}}(-tD^2)(x,g\cdot
x)]\}={\rm Tr}(T^{\lambda}_i(x)K(x)).\eqno(2.3.11)$$ \noindent{\bf
{\it Proof.}}~~Let $\{b_j\}$ be the basis of $I$, then
$\{[(\sigma(x),b_j)]\}$ is the basis of $\Gamma(S)|_U$.Then
\begin{eqnarray*}
&~&{\rm
Tr}\{(D^{\lambda}_i)_x[D{\rm{exp}}(-tD^2)(x,g\cdot x)]\}\\
&=&\sum_j\langle (D^{\lambda}_i)_x[D{\rm{exp}}(-tD^2)(x,g\cdot x)](\sigma(g\cdot x),b_j),(\sigma(g\cdot x),b_j)\rangle\\
&=&\sum_j\langle (\sigma(g\cdot x),T^{\lambda}_i(x)K(x)b_j),(\sigma(g\cdot x),b_j)\rangle\\
&=&\sum_j\langle T^{\lambda}_i(x)K(x)b_j,b_j\rangle\\
&=&{\rm Tr}(T^{\lambda}_i(x)K(x)).~~\Box
\end{eqnarray*}
\quad As in [LYZ],[Y], We define $\chi(x^\alpha
D^\beta_xe^\gamma)=|\beta|-|\alpha|+|\gamma|,~\alpha,~\beta\in {\bf
Z}^n,~\gamma\in {\bf Z_2}^n$. We also define $\chi (t)=-2$ and $\chi
(z)=1$, then we have $$\chi(zx^\alpha
D^\beta_xe^\gamma)=1+|\beta|-|\alpha|+|\gamma|.\eqno(2.3.12)$$\\
\noindent{\bf Lemma 2.12}~~{\it If $\lambda\not=0,~k\not=0$ and $k$
is even, then $\chi(t^{|\lambda|+\frac{k}{2}}T^{\lambda}_i)\leq
-2+2n'$.}\\
\noindent{\bf {\it Proof.}}~~By Lemma 3.6 in [CH], $\chi
([c(df_i)]^{[\lambda_i]})\leq 2\lambda_i$ and the simple argument in
[LYZ] shows that
$$g^*=\prod^n_{\alpha=n-2n'+1}e_{\alpha}\cdot d_1+d_2\,\eqno(2.3.13)$$
where  $\chi(d_1)\leq 0$ and $\chi(d_2)\leq 2n'-2$. So $\chi
(t^{|\lambda|+\frac{k}{2}}T^{\lambda}_i)\leq
-2(|\lambda|+\frac{k}{2})+2n'+2|\lambda|=2n'-k\leq 2n'-2.$~~~ $\Box$\\

\noindent{\bf Lemma 2.13}([Z])~~{\it Suppose $i\leq
[\frac{n}{2}]+2$. If $W$ is an odd element and $\chi(W)\leq
2i-2+2n',$ then
$$ lim_{\small t\rightarrow 0}\frac{1}{t^{\frac{3}{2}}}\left|\int_{N_\xi(\varepsilon)}
\frac{e^{-\frac{d(x',g\cdot x')^2}{4t}}}{(4\pi t)^{\frac{n}{2}}}{\rm
Tr}(W(0;x'))t^idx'\right|\leq C'\eqno(2.3.14)$$
 for some constant $C'>0$; where in the
$W(x'';x')$, $x''$ stands for tangential coordinates and $x'$
stands for normal coordinates.}\\

\noindent{\bf Theorem 2.14}~~{\it When $t\rightarrow 0$,}
$$t^{|\lambda|+\frac{k}{2}}{\rm
Tr}\{D^{\lambda}_iDe^{-tD^2}\}\sim
O(t^{\frac{1}{2}}).\eqno(2.3.15)$$ \noindent{\bf {\it Proof.}}~By
Lemma 2.11, we only need to prove that when $t\rightarrow 0$,
$$t^{|\lambda|+\frac{k}{2}}\int_{N_\xi(\varepsilon)}{\rm
Tr}(T^{\lambda}_i(0,x')K(0,x'))dx'\sim
O(t^{\frac{1}{2}}).\eqno(2.3.16)$$
 Set
$$h(x)=1+\frac{1}{2}z\sum\limits^n_{i=1}{x_ic(e_i)},\eqno(2.3.17)$$
where $(x_1,\cdots,x_n)$ is the normal coordinates under the frame
$E_1^{g\cdot x},\cdots, E_n^{g\cdot x}$ and we consider $h$ as
$h\phi$ where $\phi$ is a cut function about $(x_1,\cdots,x_n)$. By
[Z], we have
$$hc(e_i)h^{-1}=c(e_i)+(\chi=0);~~h(D^2-zD)h^{-1}=D^2+zu,\eqno(2.3.18)$$
where $\chi(u)\leq 0$,~ $u$ contains no $z$ and the equality
$$ztD{\rm{exp}}(-tD^2)(x,y)=h^{-1}(x){\rm{exp}}(-t(D^2+zu)(x,y))h(y)-{\rm{exp}}(-tD^2)(x,y)\eqno(2.3.19)$$
where
$${\rm{exp}}(-tD^2)(x,y)=\frac{e^{-\frac{d(x,y)^2}{4t}}}{(4\pi
t)^{\frac{n}{2}}}\left(\sum\limits^{[\frac{n}{2}]+2}_{i=0}U_it^i+o(t^{[\frac{n}{2}]+2})\right);\eqno(2.3.20)$$
$${\rm{exp}}(-t(D^2+zu))(x,y)=\frac{e^{-\frac{d(x,y)^2}{4t}}}{(4\pi
t)^{\frac{n}{2}}}\left(\sum\limits^{[\frac{n}{2}]+2}_{i=0}(U_i+zV_i)t^i+o(t^{[\frac{n}{2}]+2})\right),\eqno(2.3.21)$$
where $\chi (U_i)\leq 2i,~\chi (V_i)\leq 2(i-1)$ and $U_i,~V_i$
contains no $z$.\\
\noindent So we get:
$$tK(x)=\frac{e^{-\frac{d(x,g\cdot x)^2}{4t}}}{(4\pi
t)^{\frac{n}{2}}}\left[\sum\limits^{[\frac{n}{2}]+2}_{i=0}(\frac{1}{2}\sum_j
((dg-I)x)_jc(e_j))U_it^i+o(t^{[\frac{n}{2}]+2}) \right.$$
$$+\left.\left(\sum\limits^{[\frac{n}{2}]+2}_{i=0}W_it^i+o(t^{[\frac{n}{2}]+2})\right)\right],\eqno(2.3.22)$$
\noindent where $W_i$ is an odd element and $\chi(W_i)\leq 2(i-1)$.
Because we integrate about the orthogonal coordinates, so in (2.3.22) we use the orthogonal coordinates again.\\
\indent i)~If $k=0$, then $D^{\lambda}_i=\widetilde{dg}f^0$ and
(2.3.15) is the result in [Z]. If $\lambda=0$,
$$t^{\frac{k}{2}}D^{\lambda}_i=t^{\frac{k}{2}}[c(df^{i+1})]
\cdots [c(df^k)]\widetilde{dg}f^0[c(df^{1})]\cdots [c(df^i)].$$ By
(2.3.13), we have
$$t^{\frac{k}{2}}T^{\lambda}_i=A_1\prod^n_{\alpha=n-2n'+1}e_{\alpha}+A_2,\eqno(2.3.23)$$
where $A_1\in (\chi\leq 0)$ and $A_2\in (\chi\leq  2n'-2).$
So by (2.3.13), (2.3.22), (2.3.23) and Lemma 2.13, we get (2.3.15).\\
\indent ii)~we let $\lambda\not= 0$ and $k\not= 0$. By (2.3.22) and
Lemma 2.12, we have
$$t^{|\lambda|+\frac{k}{2}+1}T^{\lambda}_i(x)K(x)
=\frac{e^{-\frac{d(x,g\cdot x)^2}{4t}}}{(4\pi t)^{\frac{n}{2}}}
\left[\sum\limits^{[\frac{n}{2}]+2}_{i=0}\widetilde{W_i}t^i+o(t^{[\frac{n}{2}]+2})\right],
\eqno(2.3.24)$$ \noindent where $\widetilde{W_i}$ is an odd element
and $\chi(\widetilde{W_i})\leq 2i-2+2n'$. Also by Lemma 2.13
 we prove this Theorem.
~~$\Box $\\

\noindent {\bf 2.4 The Convergence of The Total Eta Invariant}\\

 \indent Let $C^1(N)$ be Banach algebra of once differentiable
 function on $N$ with the norm
 $$||f||_1:={\rm sup}_{x\in N}|f(x)|+{\rm sup}_{x\in N}||df(x)||.$$
 Let $$\phi^G=\{\phi^G_0,\cdots,\phi^G_{2q},\cdots\}$$
be an equivariant even cochains sequence in the bar complex of
$C^1(N),$ then
$$||\phi^G_{2q}||={\rm sup}_{||f_i||\leq 1;~0\leq i\leq
 2q}\{||\phi^G_{2q}(f_0,\cdots,f_{2q})||_{C(G)}\}.$$
 \noindent {\bf Definition 2.15 }~ The radius of convergence of
 $\phi^G$ is defined to be that of the power series $\sum
 q!||\phi^G_{2q}||z^q.$ The space of cochains sequence with radius
 of convergence at least $r>0$ is denoted by $C^{{\rm
 even},G}_r(C^1(N))$ (similarly define $C^{{\rm
 odd},G}_r(C^1(N))$).\\
 \noindent In general, the sequence
$$\eta^G(D)=\{\cdots,\eta^G_{2q}(D),\eta^G_{2q+2}(D),\cdots\}$$
which called total eta invariant is not an entire cochain.\\

\noindent {\bf Proposition 2.16} {\it Suppose that $D$ is invertible
with $\lambda$ the smallest positive eigenvalue of $|D|$. Then the
equivariant total eta invariant $\eta^G(D)$ has radius of
convergence $r$ satisfying the inequality: $r\geq 4\lambda^2>0$ i.e.
$\eta^G(D)\in C^{{\rm
 even},G}_{4\lambda^2}(C^1(N)).$}\\

\noindent{\bf {\it Proof.}} We can assume $q>(n+2)/2$ for only
considering convergence, where $n={\rm dim}N$. For a constant
$\delta>0$, we have
$$\sum_qq!||\eta^G_{2q}(D)||z^q\leq
\sum_qq!\int^{\infty}_{1+\delta}||\widetilde{{\bf {\rm
ch}}}_{2q}^G(tD,D)||dtz^q+\sum_qq!\int^{1+\delta}_0||\widetilde{{\bf
{\rm Ch}}}_{2q}^G(tD,D)||dtz^q.$$ As the proof of the proposition
1.5 in [W], the sequence
$$\sum_qq!\int^{\infty}_{1+\delta}||\widetilde{{\bf {\rm
ch}}}_{2q}^G(tD,D)||dtz^q$$ has a radius $r\geq 4\lambda^2$ of
convergence; while similar to the proof of Theorem 2.8 (2), we have
$$\sum_qq!\int^{1+\delta}_0||\widetilde{{\bf {\rm
ch}}}_{2q}^G(tD,D)||dtz^q$$ is an entire sequence. So this
proposition is correct.~~$\Box $\\
\indent For the idempotent $p\in {\cal M}_r(C^{\infty}(N))$, its
Chern character ${\rm Ch}(p)$ in entire cyclic homology is defined
by the formula (for more details see [GS]):
$${\rm Ch}(p)={\rm Tr}(p)+\sum_{k\geq 1}\frac{(-1)^k(2k)!}{k!}{\rm
Tr}_{2k}((p-\frac{1}{2})\otimes \overline{p}^{\otimes 2k})
\eqno(2.4.1)$$ where
$${\rm Tr}_{2k}:~{\cal M}_r(C^{\infty}(N))\otimes\left({\cal
M}_r(C^{\infty}(N))/{\cal M}_r({\bf C})\right)^{\otimes
2k}\rightarrow C^{\infty}(N)\otimes(C^{\infty}(N)/{\bf C})^{\otimes
2k}$$ is the generalized trace map. Let
$$||dp||=||[D,p]||=\sum_{i,j}||dp_{i,j}||\eqno(2.4.2)$$
where $p_{i,j}~(1\leq i,j\leq r)$ is the entry of $p$. \\

\noindent {\bf Proposition 2.17} {\it We assume that
$||dp||<\lambda$, then
the pairing $\langle\eta^G(D),Ch(p)\rangle$ is well-defined.}\\

\noindent{\bf {\it Proof.}} By (2.4.2), using the same method as
Proposition 2.16, we can prove this
proposition.~~$\Box $\\

\section {The Proof of (1.1)}

\quad In this section, we will prove (1.1) by using the method in
[G2]. Let
$$C^{\infty}_G(N)=\{f\in C^{\infty}(N)|f(g.x)=f(x),~{\rm for~
any}~g\in G~ {\rm and}~ x\in N\}.$$
 \noindent Suppose that $D$ is
invertible with $\lambda$ the smallest eigenvalue of $|D|$, and
$p=p^*=p^2\in {\cal M }_r(C^{\infty}_G(N))$ is a idempotent which
satisfying $||dp||<\lambda$. Let
$$p(D\otimes I_r)p:~p(H\otimes {\bf C^r})=L^2(N,S\otimes
p({\bf C^r}))\rightarrow L^2(N,S\otimes p({\bf C^r}))$$ be the Dirac
operator with coefficients from $F=p({\bf C^r})$. We denote
$$\widetilde{dg}\otimes I_r:~L^2(N,S\otimes
p({\bf C^r}))\rightarrow L^2(N,S\otimes p({\bf C^r}))$$ \noindent
still by $g$. Since $p\in {\cal M }_r(C^{\infty}_G(N))$,
 we have $$g[p(D\otimes I_r)p]=[p(D\otimes I_r)p]g.$$
{\bf Theorem 3.1} {\it Under the assumption as above, we have
 $$\frac{1}{2}\eta^G(p(D\otimes
 I_r)p)=\langle\eta^G(D),{\rm Ch}(p)\rangle,\eqno(3.1)$$
where the left term is the equivariant Atiyah-Patodi-Singer eta
invariant.}\\
 \indent Let
 $${\bf D}=\left[\begin{array}{lcr}
  \ 0 &  -D\otimes I_r \\
    \  D\otimes I_r  & 0
\end{array}\right];~
p=\left[\begin{array}{lcr}
  \ p &  0 \\
    \ 0  & p
\end{array}\right];~
\sigma=i\left[\begin{array}{lcr}
  \ 0 &  I_r \\
    \  I_r  & 0
\end{array}\right];~
g=\left[\begin{array}{lcr}
  \  g &  0 \\
    \ 0  & g
\end{array}\right],$$
be operators from $H\otimes {\bf C^r}\oplus H\otimes {\bf C^r}$ to
itself, then ${\bf D}$ is skew-adjoint and
$${\bf D}\sigma=-\sigma {\bf D};~~\sigma p=p\sigma.$$
Moreover ${\bf D}e^{t{\bf D}^2}$ and $e^{t{\bf D}^2}~ (t>0)$ are
traceclass, so $(C^{\infty}_G(N),H\otimes {\bf C^r}\oplus H\otimes
{\bf C^r},{\bf D})$ is a degree-1 Fredholm module in the sense of
[G2] (see Definition 2.1 in [G2]). For $u\in [0,1]$, Let
$$D_u=(1-u)D+u[pDp+(1-p)D(1-p)]=D+u(2p-1)[D,p],$$
then
$${\bf
D}_u=\left[\begin{array}{lcr}
  \  0 &  -D_u \\
    \ D_u  & 0
\end{array}\right]
={\bf D}+u(2p-1)[{\bf D},p].$$ \noindent We consider a family of
Fredholm modules on
 $[0,1]\times {\bf R}\times [0,\infty)$, parameterized by
$(u,s,t)$,
$$\widetilde{{\bf D}}=t^{\frac{1}{2}}{\bf
D}_u+s\sigma(p-\frac{1}{2}),$$ then $\widetilde{{\bf
D}}^*=-\widetilde{{\bf D}}$ and $g\widetilde{{\bf
D}}=\widetilde{{\bf D}}g.$ Let $A=d+\widetilde{{\bf D}}$ be a
superconnection on the trivial infinite dimensional superbundle with
base $[0,1]\times {\bf R}\times [0,\infty)$ and fibre $H\otimes {\bf
C^r}\oplus H\otimes {\bf C^r}.$ By [G2], we have
$$(d+\widetilde{{\bf D}})^2=t{\bf
D}_u^2-s^2/4-(1-u)t^{\frac{1}{2}}s\sigma[{\bf D},p]+ds\sigma
(p-\frac{1}{2})+t^{\frac{1}{2}}du(2p-1)[{\bf
D},p]+\frac{1}{2}t^{-\frac{1}{2}}dt{\bf D}_u.\eqno(3.2)$$ \noindent
We also consider $A$ as $A_t$, which is a family superconnection
parameterized  by $t$ on trivial superbundle with base $[0,1]\times
{\bf R}$ and fibre $H\otimes {\bf C^r}\oplus H\otimes {\bf C^r}.$ By
Duhamel principle and $gA_t=A_tg$, then
\begin{eqnarray*}
 \int_0^{+\infty}{\rm Str}(ge^{A^2})&=&\int_0^{+\infty}\int_0^1{\rm Str}
\left[ge^{s_0A_t^2}\frac{1}{2}t^{-\frac{1}{2}}dt{\bf
D}_ue^{(1-s_0)A^2_t}\right]ds_0\\
&=&\int_0^{+\infty}{\rm
Str}\left[g\frac{dA_t}{dt}e^{A_t^2}\right]dt.
\end{eqnarray*}
By $gA_t=A_tg$, similar to Theorem 9.23 in [BGV], we have,
\begin{eqnarray*}
d\int_0^{+\infty}{\rm Str}(ge^{A^2})&=&d\int_0^{+\infty}{\rm
Str}\left[g\frac{dA_t}{dt}e^{A_t^2}\right]dt\\
&=&{\rm lim}_{t\rightarrow +\infty}{\rm Ch}_g(A_t)-{\rm
lim}_{t\rightarrow 0}{\rm Ch}_g(A_t).
\end{eqnarray*}
Using Duhamel principle and similar to the proof of Lemma 1.1 in
[Wu], we have
$${\rm
lim}_{t\rightarrow +\infty}{\rm Ch}_g(A_t)=0.\eqno(3.3)$$ Using
Duhamel principle, as the proof of Section 2.3, we get
$${\rm
lim}_{t\rightarrow 0}{\rm Ch}_g(A_t)=0.\eqno(3.4)$$
 \noindent Let $\Gamma_u=\{u\}\times {\bf R}\subset [0,1]\times {\bf
 R}$ be a contour oriented in the direction of increasing $s$
 and $\gamma_s=[0,1]\times \{s\}$
 be a contour oriented in the direction of increasing $u$ . By
 Stokes theorem, then
 $$0=\int_{[0,1]\times {\bf
 R}}d\int_0^{+\infty}{\rm
 Str}(ge^{A^2})=\left(\int_{\Gamma_1}-\int_{\Gamma_0}
-\int_{\gamma_{+\infty}}+\int_{\gamma_{-\infty}}\right)
\left[\int_0^{+\infty}{\rm
 Str}(ge^{A^2})\right].$$
As the proof of (3.3), we have for some constant $C>0$,
$$\int_{\gamma_s}\int_0^{+\infty}{\rm
 Str}(ge^{A^2})\sim O(e^{-cs^2}).$$
So $$\int_{\Gamma_0}\int_0^{+\infty}{\rm
 Str}(ge^{A^2})=\int_{\Gamma_1}\int_0^{+\infty}{\rm
 Str}(ge^{A^2}).$$
By Duhamel principle and (3.2), when $u=0$, we have
$$\int_{\Gamma_0}\int_0^{+\infty}{\rm
 Str}(ge^{A^2})=\sum^{\infty}_{k=0}\int_{-\infty}^{+\infty}
\int_0^{+\infty}e^{-s^2/4}~~~~~~~~~~~~~~~~~~~~~~~$$
$$\times\int_{\triangle_k}{\rm
 Str}\left\{e^{t_0t{\bf D}^2}\left[-t^{\frac{1}{2}}s\sigma[{\bf D},p]
+ds\sigma (p-\frac{1}{2})+(d\sqrt{t}){\bf D}\right]\right.$$
 $$~~~~~~~~~~~~~~~\times\left. e^{t_1t{\bf D}^2}
\cdots\left[-t^{\frac{1}{2}}s\sigma[{\bf D},p] +ds\sigma
(p-\frac{1}{2})+(d\sqrt{t}){\bf D}\right] e^{t_kt{\bf
D}^2}g\right\}dt_0\cdots dt_k.\eqno(3.5)$$ \noindent Expanding (3.5)
in powers of $s$, we will integrate about $dsdt$ and
$$\int_{-\infty}^{+\infty}e^{-s^2/4}s^{2k+1}ds=0,$$
so we only keep terms with one factor of $ds$, one factor of
$d\sqrt{t}$, and odd number of factors of $\sigma$, using (2.2), we
get (3.5) equals
$$\sum^{\infty}_{l=0}\int_{-\infty}^{+\infty}e^{-s^2/4}s^{2l}ds
\int_0^{+\infty}t^ld\sqrt{t}
\left\{\sum_{i=0}^{2l}\sum_{j=0}^{2l-i}\times\right.$$
$$\left[-\langle
1,
\begin{array}{c}
\ \\
\ \underbrace{\sigma[{\bf D},p],\cdots,\sigma[{\bf D},p]}
\\
\ i~{\rm times}
\end{array}
,\sigma(p-\frac{1}{2}),
\begin{array}{c}
\ \\
\ \underbrace{\sigma[{\bf D},p],\cdots,\sigma[{\bf D},p]}
\\
\ j~{\rm times}
\end{array}
,{\bf D},\right.$$ $$
\begin{array}{c}
\ \\
\ \underbrace{\sigma[{\bf D},p],\cdots,\sigma[{\bf D},p]}
\\
\ (2l-i-j)~{\rm times}
\end{array}
\rangle_{\sqrt{t}}(g)+\langle 1,
\begin{array}{c}
\ \\
\ \underbrace{\sigma[{\bf D},p],\cdots,\sigma[{\bf D},p]}
\\
\ i~{\rm times}
\end{array}
,$$
$$\left.\left.{\bf D},
\begin{array}{c}
\ \\
\ \underbrace{\sigma[{\bf D},p],\cdots,\sigma[{\bf D},p]}
\\
\ j~{\rm times}
\end{array}
,\sigma(p-\frac{1}{2}),
\begin{array}{c}
\ \\
\ \underbrace{\sigma[{\bf D},p],\cdots,\sigma[{\bf D},p]}
\\
\ (2l-i-j)~{\rm times}
\end{array}
\rangle_{\sqrt{t}}(g)\right]\right\}.\eqno(3.6)$$ The following
equality is the equivariant case of Lemma 2.2 (2) in [GS]. Assume
$gD=Dg$ and $gA_i=A_ig$ for $0\leq i\leq n$, then
$$
\langle
A_0,\cdots,A_n\rangle_D(g)=\sum_{i=0}^n(-1)^{(|A_0|+\cdots+|A_i|)
(|A_{i+1}+\cdots+|A_n|)}
\langle1,A_{i+1},\cdots,A_n,A_0,\cdots,A_i\rangle_D(g).\eqno(3.7)$$
By (3.7), Lemma 3.2 in [G2] and $ {\bf D}=\left[\begin{array}{lcr}
  \  0 &  -D \\
    \ D & 0
\end{array}\right]$, then (3.6) equals\\
$$-\sum^{\infty}_{l=0}\int_{-\infty}^{+\infty}e^{-s^2/4}s^{2l}ds
\int_0^{+\infty}t^l\sum_{j=0}^{2l}\times~~~~~~~~~~~~~~~~~~~~~~~~~$$
 $$\langle
\sigma(p-\frac{1}{2}), \begin{array}{c}
\ \\
\ \underbrace{\sigma[{\bf D},p],\cdots,\sigma[{\bf D},p]}
\\
\ j~{\rm times}
\end{array}
,{\bf D},
\begin{array}{c}
\ \\
\ \underbrace{\sigma[{\bf D},p],\cdots,\sigma[{\bf D},p]}
\\
\ (2l-j)~{\rm times}
\end{array}
\rangle_{\sqrt{t}}(g)d\sqrt{t}$$
\begin{eqnarray*}
&=&-\sum^{\infty}_{l=0}\frac{2l!}{l!}\sqrt{4\pi} \int_0^{+\infty}
\sum_{j=0}^{2l}\\
& &~~~~~~~~~\times\langle \sigma(p-\frac{1}{2}),
\begin{array}{c}
\ \\
\ \underbrace{\sigma[t{\bf D},p],\cdots,\sigma[t{\bf D},p]}
\\
\ j~{\rm times}
\end{array}
,{\bf D},
\begin{array}{c}
\ \\
\ \underbrace{\sigma[t{\bf D},p],\cdots,\sigma[t{\bf D},p]}
\\
\ (2l-j)~{\rm times}
\end{array}
\rangle_{t}(g)dt\\
&=&-\sum^{\infty}_{l=0}\frac{2l!}{l!}\sqrt{4\pi} \int_0^{+\infty}
\sum_{j=0}^{2l}(-1)^j\\
& &~~~~~~~~~\times \langle \sigma(p-\frac{1}{2}),
\begin{array}{c}
\ \\
\ \underbrace{[t{\bf D},p],\cdots,[t{\bf D},p]}
\\
\ j~{\rm times}
\end{array}
,{\bf D},
\begin{array}{c}
\ \\
\ \underbrace{[t{\bf D},p],\cdots,[t{\bf D},p]}
\\
\ (2l-j)~{\rm times}
\end{array}
\rangle_{t}(g)dt\\
&=&-\sum^{\infty}_{l=0}(-1)^l\frac{2l!}{l!}\sqrt{4\pi}\times
2\sqrt{-1}
\int_0^{+\infty} \sum_{j=0}^{2l}(-1)^j \\
& &~~~~~~~~~~~~\times\langle (p-\frac{1}{2}),
\begin{array}{c}
\ \\
\ \underbrace{[tD,p],\cdots,[tD,p]}
\\
\ j~{\rm times}
\end{array}
,D,
\begin{array}{c}
\ \\
\ \underbrace{[tD,p],\cdots,[tD,p]}
\\
\ (2l-j)~{\rm times}
\end{array}
\rangle_{t}(g)dt\\
&=&-4\sqrt{-1}\pi[\langle \eta^G(D),{\rm
Ch}(p)\rangle(g)-\frac{1}{2}\langle \eta^G(D),{\rm rk}(p){\rm
Ch}_*(1)\rangle(g)].~~~~~~(3.8)
\end{eqnarray*}
When $u=1$, using Duhamel principle, then
$$\int_{\Gamma_1}\int_0^{+\infty}{\rm
 Str}(ge^{A^2})=\int_{-\infty}^{+\infty}\int_0^{+\infty}
e^{-s^2/4}\sum_{k=0}^{\infty}\int_{\triangle_k}{\rm
 Str}\left\{e^{t_0t{\bf D_1}^2}\right.$$
 $$\times\left.\left[ds\sigma (p-\frac{1}{2})
 +d\sqrt{t}{\bf D_1}\right]
e^{t_1t{\bf D_1}^2}\cdots \left[ds\sigma
(p-\frac{1}{2})+d\sqrt{t}{\bf D_1}\right] e^{t_kt{\bf
D_1}^2}g\right\}dt_0\cdots dt_k.\eqno(3.9)$$ \noindent By the same
reason as $u=0$, we have $k=2$ and (3.9) equals
$$\int_{-\infty}^{+\infty}\int_0^{+\infty}
e^{-s^2/4}\left\{\int_{\triangle_2 }{\rm
 Str}\left[e^{t_0t{\bf D_1}^2}ds\sigma (p-\frac{1}{2})
e^{t_1t{\bf D_1}^2}d\sqrt{t}{\bf D_1}e^{t_2t{\bf
D_1}^2}g\right]dt_0dt_1 dt_2\right.$$
$$\left.+\int_{\triangle_2} {\rm
 Str}\left[e^{t_0t{\bf D_1}^2}d\sqrt{t}{\bf D_1}
e^{t_1t{\bf D_1}^2}ds\sigma (p-\frac{1}{2})e^{t_2t{\bf
D_1}^2}g\right]dt_0dt_1 dt_2\right\}$$
$$=-\int_{-\infty}^{+\infty}e^{-s^2/4}ds
\int_0^{+\infty}{\rm
 Str}[\sigma (p-\frac{1}{2})
{\bf D_1}e^{t{\bf D_1}^2}g]d\sqrt{t}~~~~~~~~~~~~~~~~~~~~~~~~~$$
$$=-2i\int_{-\infty}^{+\infty}e^{-s^2/4}ds
\int_0^{+\infty}{\rm
 Tr}[(p-\frac{1}{2})
D_1e^{-t D_1^2}g]d\sqrt{t}~~~~~~~~~~~~~~~~\eqno(3.10)$$ \noindent
Let
$$D_p=p(D\otimes I_r)p:~p(H\otimes {\bf C^r})\rightarrow
p(H\otimes {\bf C^r}).$$
\par
\noindent{\bf Lemma 3.2}~ $$\int_0^{+\infty}{\rm
 Tr}[g
D_1e^{-t D_1^2}]d\sqrt{t}=\int_0^{+\infty}{\rm
 Tr}[g
De^{-t D^2}]d\sqrt{t}.\eqno(3.11)$$ \noindent {\it Proof.} Let
$A=(2p-1)dp$, using the property of trace, then
\begin{eqnarray*}
\frac{d}{du}\eta^G(D_u)(g)&=&\frac{2}{\sqrt{\pi}}
\int_0^{+\infty}\frac{d}{du}{\rm
 Tr}[g
D_ue^{-t D_u^2}]d\sqrt{t}\\
&=&\frac{2}{\sqrt{\pi}}\int_0^{+\infty} {\rm
 Tr}[g
Ae^{-t D_u^2}-2tg AD_u^2e^{-t D_u^2}]d\sqrt{t}\\
&=&\frac{2}{\sqrt{\pi}}\int_0^{+\infty}\frac{d}{dt}[t^{\frac{1}{2}}{\rm
 Tr}(g
Ae^{-t D_u^2})]dt\\
&=& \frac{2}{\sqrt{\pi}}t^{\frac{1}{2}}{\rm
 Tr}(g
Ae^{-t D_u^2})|^{+\infty}_0=0.
\end{eqnarray*}
In the last equality, considering $\chi(A)=1$, we use the similar
trick in [Z]. $\Box$\\
\noindent By (3.11),
 then (3.10) equals
\begin{eqnarray*}
&=&-2\pi i\{\frac{2}{\sqrt{\pi}} \int_0^{+\infty} {\rm
 Tr}[g
D_pe^{-t D_p^2}g]d\sqrt{t}-\frac{1}{2}\frac{2}{\sqrt{\pi}} {\rm
 Tr}[g
De^{-t D^2}]d\sqrt{t}\}\\
&=&-2\pi i[\eta^G(D_p)- \frac{1}{2}\eta^G(D\otimes I_r)]\\
&=&-2\pi i[\eta^G(D_p)- \frac{1}{2}{\rm rk}(p)\eta^G(D)]\\
&=&-2\pi i[\eta^G(D_p)-{\rm rk}(p)\langle\eta^G(D),{\rm
Ch}_{\star}(1)\rangle].~~~~~~~~~~~~~~~~~~~~(3.12)
\end{eqnarray*}
 \noindent Since (3.8) equals (3.12), then we get Theorem 3.1.\\

\section{The Equivariant Chern-Connes Character On Manifolds
With Boundary}

\quad Let $M$ be an even-dimensional compact spin manifolds with
boundary $\partial M=N $ endowed with a metric which is a product in
a collar neighborhood of $N$. Denote by $D ~(D_N)$ the Dirac
operator acting on the spinors bundle on $M~(N)$. Suppose that $G$
is a compact Lie group acting on $M$ by orientation-preserving
isometries. Let $C^{\infty}_G(M)=\{f\in C^{\infty}(M)|f$ is
independent of the normal coordinate $ x_n$ near the boundary
and $ f|_N\in C^{\infty}_G(N)\}.$\\

 \noindent {\bf Definition
4.1}~The equivariant Chern-Connes character on $M$, $\tau^G=\{
\tau^G_0,\tau^G_2,\cdots, \tau^G_{2q}\cdots \}$ is defined by
$$ \tau^G_{2q}(f^0,f^1,\cdot,f^{2q})(g):= -\eta^G_{2q}(D_N)
(f^0|_N,f^1|_N,\cdot,f^{2q}|_N)(g)+$$
$$ \frac{1}{(2q)!(2\pi\sqrt{-1})^q}\sum_{i=1}^k
\int_{F_i} \widehat{A}(TF_i) \left\{{\rm Pf}\left[2{\rm
sinh}(\Omega/{4\pi}+\frac{\sqrt{-1}{\bf
\theta}}{2})(N(F_i))\right]\right\}^{-1}f^0 df^1\wedge \cdots\wedge
df^{2q},    \eqno(4.1)$$ where $\{F_1,\cdots,F_k\}$ are components
of the fixed point set of $g$ acting on $M$. $\Omega$ is the
curvature matrix of the normal bundle $N(F_i)$ and $\theta$ is a
function matrix on $F_i$ (For details, see [LYZ]).
 $f_i\in C^{\infty}_G(M)~(0\leq i\leq 2q)$,
$\eta^G_{2q}(D_N)$ is the equivariant $\eta$-cochain defined in
Section 2 and
 $ \widehat{A}(TF_i)$ is the $\hat{A}$-polynomial of curvature $R_{F_i}$ of $F_i$.\\

 \noindent {\bf Proposition 4.2}~{\it The equivariant Chern-Connes character is
$b-B$ closed (for the definitions of $b,~B$, see [FGV]), i.e.,}
$$b\tau^G_{2q-2}+B\tau^G_{2q}=0.\eqno(4.2) $$
\noindent {\bf Proof.}~Firstly, we have the equivariant version of
the corollary 2.5 in [GS]. Let $gD_N=D_Ng$ for $g\in G$ and $f^i\in
C^{\infty}_G(N)~(0\leq i \leq 2q+1)$, we have
$$ -\frac{d}{dt}{\bf {\rm ch}}_{2q+1}^G(tD_N)(f^0,\cdots ,
f^{2q+1})(g)= b\widetilde{{\bf {\rm
ch}}_{2q}^G}(tD_N,D_N)(f^0,\cdots , f^{2q+1})(g)$$
$$+B\widetilde{{\bf {\rm ch}}_{2q+2}^G}(tD_N,D_N)(f^0,\cdots ,
f^{2q+1})(g). \eqno(4.3)$$
 \noindent Similar to the discussion in [CM], we have
$$ {\rm lim}_{t\rightarrow \infty}{\bf {\rm ch}}_{2q+1}^G(tD_N)(f^0,\cdots ,
f^{2q+1})(g)=0. \eqno(4.4)$$ \noindent By (4.3) and (4.4), then
$$\frac{1}{\Gamma(\frac{1}{2})}{\rm
lim}_{t\rightarrow 0} {\bf {\rm ch}}_{2q+1}^G(tD_N)(f^0,\cdots ,
f^{2q+1})(g)=[b\eta^G_{2q}(D_N)+B\eta^G_{2q+2}(D_N)](f^0,\cdots ,
f^{2q+1})(g). \eqno(4.5)$$ \noindent Similar to the computation in
[CH], we get
$$\frac{1}{\Gamma(\frac{1}{2})}{\rm
lim}_{t\rightarrow 0} {\bf {\rm ch}}_{2q+1}^G(tD_N)(f^0,\cdots ,
f^{2q+1})(g)= \frac{1}{(2q+1)!(2\pi\sqrt{-1})^{q+1}}\sum_{i=1}^k
\int_{\partial F_i} \widehat{A}(T\partial F_i)$$
$$\times \left\{{\rm
Pf}\left[2{\rm sinh}(\Omega/{4\pi}+\frac{\sqrt{-1}{\bf
\theta}}{2})(N(\partial F_i))\right]\right\}^{-1}f^0d_Nf^1\wedge
\cdots\wedge d_Nf^{2q+1}.\eqno(4.6)$$ \noindent By (4.5) and (4.6),
similar to
the discussion of [Wu], we get (4.2).~~~~~~$\Box$\\

 \noindent {\bf Remark}~ By $f\in C^{\infty}_G(M)$, so $b^G~(B^G)$ defined
in [KL] is $b~(B)$.\\

\indent Let $C^1_G(M)$ be the Banach algebra, the completion of
$C^{\infty}_G(M)$ under the norm $||.||_1$ defined in Section 2. Let
$\phi^G=\{ \phi_0^G ,\cdots,\phi_{2q}^G ,\cdots,\}$ be the
 an equivariant even cochains sequence in the bar complex of $C^1_G(M)$. We call that
 $\phi^G$ have radius of convergence at least $r>0$ relative to $N$ if $\phi^G$ can be
written as a sum of two cochains $\phi^G=\phi^{(1),G }+\phi^{(2),G}
$ with $\phi^{(1),G } $ entire and $\phi^{(2),G} $ supported on $N$
such that $\phi^{(2),G}\in C^{\star,G}_r(C^1_G(N)).$  Then we have a
corollary of Proposition
2.6 and Proposition 2.7. \\

 \noindent {\bf Proposition 4.3}~{\it The equivariant Chern-Connes character $\tau^G$
 has radius of convergence  at least $4\lambda^2$, where $\lambda$ is the smallest positive
eigenvalue of the invertible operator $D_N$. For a selfadjoint
idempotent $p\in M_r(C^{\infty}_G(M))$ such that
$||d(p|_N)||<\lambda$, then
the pairing $\langle\tau^G,{\rm Ch}(p)\rangle(g)$ is well-defined.}\\

\section{The Index Pairing}

\quad In this section, we will give the main result of this
paper.\\
\indent Let $M$ be a smooth connected compact manifold with smooth
 compact boundary $N$. Assume $M$  has even dimension, is oriented
 and spin. Let
 $$C_1(N)=N\times (0,1];~~Z=M\cup_{N\times \{1\}}C_1(N),$$
and $\cal{U}$ be a collar neighborhood of $N$ in $M$. For
$\varepsilon>0$, we take a metric $g^{\varepsilon}$ of $Z$ such that
on ${\cal{U}}\cup_{N\times \{1\}}C_1(N)$
$$g^{\varepsilon}=\frac{dr^2}{\varepsilon}+r^2g^{N}.$$
\noindent Let $S=S^+\oplus S^-$ be spinors bundle associated to
$(Z,g^{\varepsilon})$ and $H^{\infty}$ be the set
$\{\xi\in\Gamma(Z,S)|~\xi ~{\rm and~ its~ derivatives~ are~ zero~
near~ the~ vertex~ of~ cone~}\}.$ Denote by $L^2_c(Z,S)$ the
$L^2-$completion of $H^{\infty}$ (similar define $L^2_{c}(Z,S^+)$
and $L^2_{c}(Z,S^-)$). Let
$$D_{\varepsilon}:~H^{\infty}\rightarrow H^{\infty};
~~D_{+,\varepsilon}:~H_+^{\infty}\rightarrow H_-^{\infty},$$ be the
Dirac operators associated to $(Z,g^{\varepsilon})$ which are
Fredholm operators for the sufficient small $\varepsilon.$ Suppose
that $G$ is a compact connected Lie
group acting on $M$ by orientation-preserving isometries.\\

\noindent ${\bf (H_1)}$~~ Assume that the boundary Dirac operator
$D_N$ is invertible and $p=p^*=p^2\in {\cal
M}_r(C^{\infty}_G(M))\subset {\cal M}_r( C^{\infty}(Z))$ such that
$||d(p|_N)||<\lambda,$ where $\lambda$ is the smallest positive
eigenvalue of $|D_N|$.\\

Consider
$$D_{p,\varepsilon}^+:=p(D_{\varepsilon}^+\otimes I_r)p:~p(
L^2_{c}(Z,S^+)\otimes {\bf C^r})\rightarrow p(L^2_{c}(Z,S^-)\otimes
{\bf C^r}),$$ which is the Dirac operator with the coefficient from
$G$-vector bundle $p({\bf C^r})$ over
$Z$. We also assume that \\

\noindent ${\bf (H_2)}$~~ For any $g\in G$, there are lifts of $g$:
$$g_1:~L^2(N,S_N\otimes {\rm Im}(p|_N))\rightarrow L^2(N,S_N\otimes {\rm
Im}(p|_N));$$
$$g_2:~L^2_{c}(Z,S\otimes {\rm Im}(p))\rightarrow
L^2_{c}(Z,S\otimes {\rm Im}(p)),$$ \noindent which commute with
$D_N$ and $D_{p,\varepsilon}$ respectively.\\

 \noindent Under the assumption  ${\bf (H_2)}$, we define
 $${\rm Ind}_gD^+_{p,\varepsilon}={\rm Tr}g|_{{\rm
 ker}D^+_{p,\varepsilon}}-{\rm Tr}g|_{{\rm
 ker}D^-_{p,\varepsilon}}.$$
 \noindent Similar to the discussion of [W, p.165], by ${\bf (H_1)}$
then
$$D_{N,p|_N}=p|_N(D_N\otimes I_r)p|_N:~ L^2(N,S_N\otimes {\rm
Im}(p|_N))\rightarrow L^2(N,S_N\otimes {\rm Im}(p|_N))$$
 \noindent is invertible. Let ${\rm dim}M=2m$. If we take the
connection $pd$ of the
 bundle ${\rm Im}(p)$, by $g={\rm id}$ on
 ${\rm Im}p|_{N_q}$, we get (see [FGV])
$${\rm Ch}_g({\rm
Im}(p))=\sum_{k=0}^{\infty}(-\frac{1}{2\pi\sqrt{-1}})^k\frac{1}{k!}{\rm
Tr}[p(dp)^{2k}].\eqno (5.1)$$
 \noindent So by Theorem 3.3 in [Z]
(also see [D]), Theorem 3.1 and (5.1), we get\\

\noindent {\bf Theorem 5.1}~{\it Under the assumption ${\bf (H_1)}$
and ${\bf (H_2)}$, then
$${\rm Ind}_gD^+_{p,\varepsilon}=
\sum_{r=0}^m \sum_{i=1}^k\frac{(-1)^r}{r!(2\pi\sqrt{-1})^r}
\int_{F_i}
\widehat{A}(TF_i)\times~~~~~~~~~~~~~~~~~~~~~~~~~~~~~~~~~~~~~~~~~~$$
$$\left\{{\rm Pf}\left[2{\rm
sinh}(\Omega/{4\pi}+\frac{\sqrt{-1}{\bf
\theta}}{2})(N(F_i))\right]\right\}^{-1}{\rm
Tr}[p(dp)^{2r}]-\langle\eta^G(D_N)(g),{\rm
Ch}(p)\rangle.\eqno(5.2)$$}

 \indent Let
$$ \hat{\tau}^G_{2q}(g)(f^0,f^1,\cdot,f^{2q})(g):=
\frac{1}{(2q)!(2\pi\sqrt{-1})^q}\sum_{i=1}^k \int_{F_i}
\widehat{A}(TF_i)\times~~~~~~~~~~~~$$ $$ \left\{{\rm Pf}\left[2{\rm
sinh}(\Omega/{4\pi}+\frac{\sqrt{-1}{\bf
\theta}}{2})(N(F_i))\right]\right\}^{-1}f^0\wedge df^1\wedge
\cdots\wedge df^{2q},    \eqno(5.3)$$ \noindent then by the Stokes
theorem, we have\\

$\langle\hat{\tau}^G,{\rm Ch}(p)\rangle(g)$\\

\noindent $=\langle\hat{\tau}_0^G,{\rm tr}(p)\rangle(g)+\sum_{q\geq
1}\langle\hat{\tau}_{2q}^G,\frac{(-1)^q(2q)!}{q!}{\rm
Tr}((p-\frac{1}{2})\otimes \overline{p}^{\otimes 2q})\rangle(g)$\\
$$=\sum_{q\geq 0}\frac{(-1)^q}{{q}!(2\pi\sqrt{-1})^q}\sum_{i=1}^k
\int_{F_i} \widehat{A}(TF_i)\times \left\{{\rm Pf}\left[2{\rm
sinh}(\Omega/{4\pi}+\frac{\sqrt{-1}{\bf
\theta}}{2})(N(F_i))\right]\right\}^{-1}{\rm Tr}[p(dp)^{2q}]$$
$$
 -\sum_{q\geq
1}\frac{(-1)^q}{{q}!(2\pi\sqrt{-1})^q}\sum_{i=1}^k \int_{F_i}
\widehat{A}(TF_i) \left\{{\rm Pf}\left[2{\rm
sinh}(\Omega/{4\pi}+\frac{\sqrt{-1}{\bf
\theta}}{2})(N(F_i))\right]\right\}^{-1}{\rm
Tr}[\frac{1}{2}(dp)^{2q}]$$
$$=\sum_{q\geq 0}\frac{(-1)^q}{{q}!(2\pi\sqrt{-1})^q}\sum_{i=1}^k
\int_{F_i} \widehat{A}(TF_i)\left\{{\rm Pf}\left[2{\rm
sinh}(\Omega/{4\pi}+\frac{\sqrt{-1}{\bf
\theta}}{2})(N(F_i))\right]\right\}^{-1}{\rm Tr}[p(dp)^{2q}]$$
$$ -\frac{1}{2}\sum_{q\geq
1}\frac{(-1)^q}{{q}!(2\pi\sqrt{-1})^q}\sum_{i=1}^k \int_{F_i}
d\{\widehat{A}(TF_i)\left\{{\rm Pf}\left[2{\rm
sinh}(\Omega/{4\pi}+\frac{\sqrt{-1}{\bf
\theta}}{2})(N(F_i))\right]\right\}^{-1}{\rm Tr}[p(dp)^{2q-1}]\}$$
 $$=\sum_{q\geq
0}\frac{(-1)^q}{{q}!(2\pi\sqrt{-1})^q}\sum_{i=1}^k \int_{F_i}
\widehat{A}(TF_i) \left\{{\rm Pf}\left[2{\rm
sinh}(\Omega/{4\pi}+\frac{\sqrt{-1}{\bf
\theta}}{2})(N(F_i))\right]\right\}^{-1}{\rm Tr}[p(dp)^{2q}]$$
$$ -\frac{1}{2}\sum_{q\geq
1}\frac{(-1)^q}{{q}!(2\pi\sqrt{-1})^q}\sum_{i=1}^k \int_{\partial
F_i} \widehat{A}(T\partial F_i)\times \left\{{\rm Pf}\left[2{\rm
sinh}(\Omega/{4\pi}+\frac{\sqrt{-1}{\bf \theta}}{2})(N(\partial
F_i))\right]\right\}^{-1}{\rm Tr}[p(dp)^{2q-1}].$$ \noindent So,
suppose that $g$ acting on $N$ has no fixed points,
then by Theorem 5.1 and (4.1), we have\\

\noindent {\bf Theorem 5.2}~{\it Suppose that $g$ acting on $N$ has
no fixed points. Under the assumption ${\bf (H_1)}$ and ${\bf
(H_2)}$, then}
$${\rm Ind}_gD^+_{p,\varepsilon}=\langle\tau^G(D),{\rm
Ch}(p)\rangle(g).\eqno(5.3)$$

 \noindent {\bf Remark:} Theorem 5.1 and 5.2 are
easily to extend to the case of twisting a $G$-vector bundle.\\

\noindent{\bf Acknowledgment.}~~ The author is indebted to Professor
Weiping Zhang for his guidance and very helpful discussions. He
thanks Professor Huitao Feng for his generous help and discussions.
He also thanks the referee for his careful reading
and helpful comments.\\

\noindent{\bf Reference}\\

 \noindent [APS] M. F. Atiyah, V. K. Patodi and
I. M. Singer, {\it Spectral asymmetry and Riemannian geometry},
Math. Proc. Cambridge Philos. Soc. 77 (1975), 43-69; 78 (1975),
405-432; 79 (1976), 71-99.\\
\noindent [BGV] N. Berline, E. Getzler, and M. Vergne, {\it Heat
kernels and Dirac operators}, Spring-Verlag, Berline Heidelberg,
1992.

\noindent [BF] J. M. Bismut and D. S. Freed, {\it The analysis of
elliptic families II}, Commun. Math. Phys. 107 (1986), 103-163.

\noindent [CH] S. Chern and X. Hu, {\it Equivariant Chern character
for the invariant Dirac operators}, Michigan Math. J. 44 (1997),
451-473.

\noindent [CM] A. Connes and H. Moscovici,{\it Transgression and
Chern character of finite dimensional K-cycles},
 Commun. Math. Phys. 155 (1993), 103-122.

\noindent [D] H. Donnelly, {\it Eta invariants for G-space}, Indiana
Univ. Math. J. 27 (1978), 889-918.

\noindent [F] H. Feng, {\it A note on the noncommutative Chern
character (in Chinese)}, Acta Math. Sinica 46 (2003), 57-64.

\noindent [FGV] H. Figueroa, J. Gracia-Bond\'{i}a and J.
V\'{a}rilly, {\it Elements of noncommutative geometry},
Birkh\"{a}user Boston, 2001.

\noindent [G1]E. Getzler, {\it Pseudodifferential operators on
supermanifolds and the Atiyah-Singer index theorem}, Commun. Math.
Phys. 92 (1983), 163-178.

\noindent [G2] E. Getzler, {\it Cyclic homology and the
Atiyah-Patodi-Singer index theorem}, Contemp. Math. 148 (1993),
19-45.

\noindent [GS] E. Getzler, and A. Szenes, {\it On the Chern
character of theta-summable Fredholm modules}, J. Func. Anal. 84
(1989), 343-357.

\noindent [KL] S. Klimek and A. Lesniewski, {\it Chern character in
equivariant entire cyclic cohomology}, K-Theory 4 (1991), 219-226.

\noindent [LYZ] J. D. Lafferty, Y. L. Yu and W. P. Zhang, {\it A
direct geometric proof of Lefschetz fixed point formulas}, Trans.
AMS. 329 (1992), 571-583.

\noindent [S] B. Simon, {\it Trace ideals and their applications},
London Math. Soc. Lecture Note 35, Cambridge University Press, 1979.

\noindent [Wu] F. Wu, {\it The Chern-Connes character for the Dirac
operators on manifolds with boundary}, K-Theory 7 (1993), 145-174.

\noindent [Y] Y. L. Yu, {\it Local index theorem for Dirac
operator}, Acta Math. Sinica (New Series) 3 (1987), 152-169.

\noindent [Z] W. P. Zhang, {\it A note on equivariant eta
invariants}, Proc. AMS. 108 (1990), 1121-1129.\\

\end{document}